\documentclass{amsart}
\issueinfo{00}
{0} 
{0} 
{2007} 

\newtheorem{theorem}{Theorem}[section]

\newtheorem{corollary}[theorem]{Corollary}
\newtheorem{definition}[theorem]{Definition}

\newtheorem{lemma}[theorem]{Lemma}
\newtheorem{proposition}[theorem]{Proposition}
\theoremstyle{remark}
\newtheorem{remark}[theorem]{Remark}
\numberwithin{equation}{section}

\begin{document}
\title[Bi-Legendrian structures and paracontact geometry]{Bi-Legendrian structures\\ and paracontact geometry}

\author[B. Cappelletti Montano]{Beniamino Cappelletti Montano}
 \address{Dipartimento di Matematica,
Universit\`{a} degli Studi di Bari, Via Edoardo Orabona 4, 70125
Bari, Italy}
 \email{cappelletti@dm.uniba.it}



\begin{abstract}
We study the interplays between paracontact geometry and the theory
of bi-Legendrian manifolds. We interpret the bi-Legendrian
connection of a bi-Legendrian manifold $M$ as the  paracontact
connection of a canonical paracontact structure induced on $M$ and
then we discuss many consequences of that result both for
bi-Legendrian  and for paracontact manifolds, as a classification of
paracontact metric structures. Finally new classes of examples of
paracontact manifolds are presented.
\end{abstract}

\maketitle

\section{Introduction}
A bi-Legendrian manifold is by definition a contact manifold
$(M,\eta)$ foliated by two transversal Legendrian foliations
${\mathcal F}_1$ and ${\mathcal F}_2$. More generally, if $M$ is
endowed with a pair of transversal, not necessarily integrable
Legendrian distributions, we speak of an almost bi-Legendrian
manifold. The study of such structures is rather recent in
literature, being started in the 90's by the works of M. Y. Pang, P.
Libermann, N. Jayne et alt. on Legendrian foliations (\cite{jayne1},
\cite{libermann}, \cite{pang}).

In this note we study some properties of bi-Legendrian manifolds. In
particular we recognize that the theory of bi-Legendrian manifolds
is closely linked to paracontact geometry. We recall that
 paracontact manifolds are semi-Riemannian  manifolds
which can be viewed as the odd dimensional counterpart of
paracomplex manifolds (see $\S$ \ref{preparacontact} for a precise
definition). These manifolds were introduced by S. Kaneyuki in
\cite{kaneyuki1} and then studied by other authors. More recently,
there seems to be an increasing interest in paracontact geometry,
and in particular in para-Sasakian geometry, due to its links to the
more consolidated theory of para-K\"{a}hler manifolds and
 to their role in geometry and mathematical physics (cf. e.g.
\cite{alekseevski1}, \cite{cortes1}, \cite{cortes2}, \cite{olszak},
\cite{erdem}). Many progresses in that subject have been reached by
some very recent papers of D. V. Alekseevski, S. Ivanov, S. Zamkovoy
and their collaborators (\cite{alekseevski2}, \cite{zamkovoy2},
\cite{zamkovoy}). In particular in \cite{zamkovoy} a complete
arrangement of all the theory  is obtained and  it is introduced a
canonical connection, called paracontact connection, which reveals
to be very useful in the study of paracontact manifolds.

The main result of this paper is that given an almost bi-Legendrian
structure on a contact manifold $(M,\eta)$, it is induced on $M$ a
canonical paracontact metric structure, and conversely if we start
from a paracontact metric manifold $(M,\psi,\xi,\eta,g)$, one can
construct on $M$ a canonical almost bi-Legendrian structure. This
result has many consequences, arising from  the interplays between
these two geometric structures. In particular we are able to find
new properties and new examples of both bi-Legendrian and
paracontact manifolds. More precisely we provide a classification of
paracontact metric structures based on the Pang's classification of
Legendrian foliations \cite{pang} and prove that, under some natural
assumptions of integrability, the canonical paracontact connection
of a paracontact metric manifold coincides with the bi-Legendrian
connection of the induced almost bi-Legendrian structure, introduced
 in \cite{mino1}. This last result yields a vanishing phenomenon of
certain characteristic classes in a para-Sasakian manifold,
providing an obstruction to the existence of this structure.
Moreover, we find some results on the curvature of the canonical
paracontact connection of a para-Sasakian manifold and characterize
flat para-Sasakian manifolds.

\section{Preliminaries}

\subsection{Legendrian foliations}\label{prefoliations}
A \emph{contact manifold} is a $(2n+1)$-dimensional smooth manifold
$M$ which admits a $1$-form $\eta$ satisfying
$\eta\wedge\left(d\eta\right)^n\neq 0$ everywhere on $M$. It is well
known that given $\eta$ there exists a unique vector field $\xi$,
called \emph{Reeb vector field}, such that $i_{\xi}\eta=1$ and
$i_{\xi}d\eta=0$. We denote by $\mathcal D$ the $2n$-dimensional
distribution defined by $\ker\left(\eta\right)$, called the
\emph{contact distribution}. It is easy to see that the Reeb vector
field is an infinitesimal automorphism with respect to the contact
distribution and  the tangent bundle of $M$ splits as the direct sum
$TM=\mathcal D\oplus\mathbb{R}\xi$. A Riemannian metric $G$ on $M$
is an \emph{associated metric} for a contact form $\eta$ if
$G\left(X,\xi\right)=\eta\left(X\right)$ for all
$X\in\Gamma\left(TM\right)$ and there exists a tensor field $\phi$
of type $(1,1)$ on $M$ such that $\phi^2=-I+\eta\otimes\xi$ and
$d\eta\left(X,Y\right)=G\left(X,\phi Y\right)$ for all
$X,Y\in\Gamma\left(TM\right)$. From the previous conditions it
easily follows  that $\phi\xi=0$, $\eta\circ\phi=0$ and
$\phi|_{\mathcal D}$ is an isomorphism. We refer to
$\left(\phi,\xi,\eta,G\right)$ as a \emph{contact metric structure}
and to $M$ endowed with such a structure as a \emph{contact metric
manifold}. A contact metric structure such that the Levi-Civita
connection satisfies
\begin{equation}\label{condizionesasaki}
(\nabla_{X}\phi)Y=g(X,Y)\xi-\eta(Y)X
\end{equation}
is said to be a \emph{Sasakian manifold}. For more details we refer
the reader to \cite{blair0}.

Note that the condition $\eta\wedge(d\eta)^n\neq 0$ implies that the
contact distribution is never integrable. One can prove that the
maximal dimension of an integrable subbundle of $\mathcal D$ is $n$.
This motivates the following definition.

\begin{definition}
A \emph{Legendrian distribution} on a contact manifold $(M,\eta)$ is
an $n$-dimensional subbundle $L$ of the contact distribution such
that $d\eta\left(X,X'\right)=0$ for all
$X,X'\in\Gamma\left(L\right)$. When $L$ is integrable, it defines a
\emph{Legendrian foliation} of $(M,\eta)$. Equivalently, a
Legendrian foliation of $(M,\eta)$ is a foliation of $M$ whose
leaves are $n$-dimensional $C$-totally real submanifolds of
$(M,\eta)$.
\end{definition}

\begin{remark}
Note that when $L$ is involutive the condition that $d\eta(X,X')=0$
for all $X,X'\in\Gamma(L)$ is unnecessary. Indeed by the
integrability of $L$ one has
$2d\eta(X,X')=X(\eta(X'))-X'(\eta(X))-\eta([X,X'])=0$.
\end{remark}

Legendrian foliations have been extensively investigated in recent
years from various points of views. In particular M. Y. Pang
provided a classification of Legendrian foliations by means of a
bilinear symmetric form $\Pi_{\mathcal F}$ on the tangent bundle of
the foliation, defined by $\Pi_{\mathcal
F}\left(X,X'\right)=-\left({\mathcal L}_{X}{\mathcal
L}_{X'}\eta\right)\left(\xi\right)$. He called a Legendrian
foliation $\mathcal F$  \emph{non-degenerate}, \emph{degenerate} or
\emph{flat} according to the circumstance that the bilinear form
$\Pi_{\mathcal F}$ is non-degenerate, degenerate or vanishes
identically, respectively. A geometrical interpretation of this
classification is given in the following lemma.

\begin{lemma}[\cite{jayne1}]\label{classificazione}
Let $\left(M,\phi,\xi,\eta,g\right)$ be a contact metric manifold
foliated by a Legendrian foliation $\mathcal F$. Then
\begin{enumerate}
    \item[(i)] $\mathcal F$ is flat if and only if
    $\left[\xi,X\right]\in\Gamma\left(T{\mathcal F}\right)$ for all
    $X\in\Gamma\left(T{\mathcal F}\right)$,
    \item[(ii)] $\mathcal F$ is degenerate if and only if there exist
    $X\in\Gamma\left(T{\mathcal F}\right)$ such that $\left[\xi,X\right]\in\Gamma\left(T{\mathcal F}\right)$,
    \item[(iii)] $\mathcal F$ is non-degenerate if and only if
    $\left[\xi,X\right]\notin \Gamma\left(T{\mathcal F}\right)$ for all
    $X\in\Gamma\left(T{\mathcal F}\right)$.
\end{enumerate}
\end{lemma}

An interesting subclass of non-degenerate Legendrian foliations is
given by those for which $\Pi$ is positive definite. Then any such
Legendrian foliation is called a \emph{positive definite Legendrian
foliation}.

It should be remarked that analogous definitions can be given also
for Legendrian distributions. Indeed let $L$ be a Legendrian
distribution and define a bilinear map $\Pi_L$ on $L$ by setting
$\Pi_{L}(X,X')=-\left({\mathcal L}_{X}{\mathcal
L}_{X'}\eta\right)\left(\xi\right)$ for all $X,X'\in\Gamma(L)$. Then
an easy computation yields $\Pi_{L}(X,X')=-\eta([\xi,X],X')$. Next,
the condition $d\eta(X,X')=0$ implies  that
$[X,X']\in\Gamma({\mathcal D})$ and hence also
$[[X,X'],\xi]\in\Gamma({\mathcal D})$. Thus
$\Pi_{L}(X,X')=-\eta([[\xi,X],X'])=\eta([[X,X'],\xi])+\eta([[X',\xi],X])=-\eta([[\xi,X'],X])=\Pi_{L}(X',X)$
and we conclude that $\Pi_L$ is symmetric. Therefore we can speak of
non-degenerate, degenerate and flat Legendrian distributions.

By an \emph{almost bi-Legendrian manifold} we mean a contact
manifold $(M,\eta)$ endowed with two transversal Legendrian
distributions $L_1$ and $L_2$. Thus, in particular, the tangent
bundle of $M$ splits up as the direct sum $TM=L_1\oplus
L_2\oplus\mathbb{R}\xi$. When both $L_1$ and $L_2$ are integrable we
speak of a \emph{bi-Legendrian manifold} (\cite{mino2}). An (almost)
bi-Legendrian manifold is said to be flat, degenerate or
non-degenerate if and only if both the Legendrian distributions are
flat, degenerate or non-degenerate, respectively.

In \cite{mino2} to any almost bi-Legendrian manifold it has been
attached  a canonical connection which plays an important role in
the study of almost bi-Legendrian manifolds.

\begin{theorem}[\cite{mino1}]\label{biconnection}
Let $(M,\eta,L_1,L_2)$ be an almost bi-Legendrian manifold. There
exists a unique connection ${\nabla}^{bl}$ such that
\begin{enumerate}
  \item[(i)] ${\nabla}^{bl} L_1\subset L_1$, \ ${\nabla}^{bl} L_2\subset L_2$, \ ${\nabla}^{bl}\left(\mathbb{R}\xi\right)\subset\mathbb{R}\xi$,
  \item[(ii)] ${\nabla}^{bl} d\eta=0$,
  \item[(iii)] ${T}^{bl}\left(X,Y\right)=2d\eta\left(X,Y\right){\xi}$ \ for all
  $X\in\Gamma(L_1)$, $Y\in\Gamma(L_2)$,\\
${T}^{bl}\left(X,\xi\right)=[\xi,X_{L_1}]_{L_2}+[\xi,X_{L_2}]_{L_1}$
  \ for all   $X\in\Gamma\left(TM\right)$,
\end{enumerate}
where ${T}^{bl}$ denotes the torsion tensor of ${\nabla}^{bl}$ and
$X_{L_1}$ and $X_{L_2}$ the projections of $X$ onto the subbundles
$L_1$ and $L_2$ of $TM$, respectively, according to the
decomposition $TM=L_1\oplus L_2\oplus\mathbb{R}\xi$.
\end{theorem}

Such \ a \ connection \ is \ called \ the \ \emph{bi-Legendrian \
connection} \ of  \ the \ almost bi-Legendrian manifold
$(M,\eta,L_1,L_2)$. Its explicit definition is the following
(\cite{mino1}). Let $H:TM\longrightarrow TM$ be the operator defined
by setting, for all $Z,Z'\in\Gamma\left(TM\right)$,
$H\left(Z,Z'\right)$ the unique section of $\mathcal{D}$ satisfying
$i_{H\left(Z,Z'\right)}d\eta|_{\mathcal{D}}=\left({\mathcal{L}}_{Z}i_{Z'}d\eta\right)|_{\mathcal{D}}$.
Then we set $\nabla^{bl}\xi:=0$ and, for any $X\in\Gamma(L_1)$,
$Y\in\Gamma(L_2)$,  $Z\in\Gamma\left(TM\right)$,
\begin{gather*}
\nabla^{bl}_{Z}X:=H\left(Z_{L_1},X\right)_{L_1}+\left[Z_{L_2},X\right]_{L_1}+\left[Z_{\mathbb
R\xi},X\right]_{L_1},\\
\nabla^{bl}_{Z}Y:=H\left(Z_{L_2},Y\right)_{L_2}+\left[Z_{L_1},Y\right]_{L_2}+\left[Z_{\mathbb
R\xi},Y\right]_{L_2}.
\end{gather*}

Further properties of that connection are collected in the following
proposition.

\begin{proposition}[\cite{mino2}]\label{proprieta}
Let $(M,\eta,L_1,L_2)$ be an almost bi-Legendrian manifold and let
${\nabla}^{bl}$ denote the corresponding bi-Legendrian connection.
Then the $1$-form $\eta$ is ${\nabla}^{bl}$-parallel and the torsion
tensor field is given by
${T}^{bl}\left(X,X'\right)=-\left[X,X'\right]_{L_2}$ for all
$X,X'\in\Gamma(L_1)$ and
${T}^{bl}\left(Y,Y'\right)=-\left[Y,Y'\right]_{L_1}$  for all
$Y,Y'\in\Gamma(L_2)$. Moreover, if $L_1, L_2$ are integrable and
flat, the curvature tensor field of $\nabla^{bl}$ vanishes along the
leaves of the foliations defined by $L_1$, $L_1\oplus\mathbb{R}\xi$,
$L_2$ and $L_2\oplus\mathbb{R}\xi$.
\end{proposition}

\subsection{Paracontact manifolds}\label{preparacontact}
A $(2n+1)$-dimensional smooth manifold $M$ has an \emph{almost
paracontact structure} (\cite{kaneyuki1}) if it admits a
$(1,1)$-tensor field $\psi$, a vector field $\xi$ and a $1$-form
$\eta$ satisfying the following conditions
\begin{enumerate}
  \item[(i)] $\eta(\xi)=1$, \ $\psi^2=I-\eta\otimes\xi$,
  \item[(ii)] denoted by $\mathcal D$ the $2n$-dimensional distribution generated by
  $\eta$, the tensor field $\psi$ induces an almost paracomplex
  structure on each fibre on $\mathcal D$.
\end{enumerate}
Recall that an almost paracomplex structure on a $2n$-dimensional
smooth manifold is a tensor field $J$ of type $(1,1)$ such that
$J\neq I$, $J^2=I$ and the eigendistributions $T^+, T^-$
corresponding to the eigenvalues $1, -1$ of $J$, respectively, have
equal dimension $n$.

As an immediate consequence of the definition one has that
$\psi\xi=0$, $\eta\circ\psi=0$ and the field of endomorphisms $\psi$
has constant rank $2n$. Any almost paracontact manifold admits a
semi-Riemannian metric $g$ such that
\begin{equation}\label{compatibile}
g(\psi X,\psi Y)=-g(X,Y)+\eta(X)\eta(Y)
\end{equation}
for all  $X,Y\in\Gamma(TM)$. Then $(M,\psi,\xi,\eta,g)$ is called an
\emph{almost paracontact metric manifold}. Note that any such
semi-Riemannian metric is necessarily of signature $(n+1,n)$. If in
addition $d\eta(X,Y)=g(X,\psi Y)$ for all $X,Y\in\Gamma(TM)$,
$(M,\psi,\xi,\eta,g)$ is said to be a \emph{paracontact metric
manifold}.

On an almost paracontact manifold one defines the tensor field
$N^{(1)}=N_{\psi}-d\eta\otimes\xi$, where $N_{\psi}$ is the
Nijenhuis tensor of $\psi$, defined as
$N_{\psi}(X,Y)=\psi^2[X,Y]+[\psi X,\psi Y]-\psi[\psi
X,Y]-\psi[X,\psi Y]$. If $N^{(1)}$ vanishes identically the almost
paracontact manifold  is said to be \emph{normal}. We prove the
following characterization of the normality in terms of foliations.

\begin{proposition}\label{parasasakian}
Let $(M,\psi,\xi,\eta)$ be an almost paracontact manifold. Let $T^+$
and $T^-$ be the eigendistributions of $\psi|_{\mathcal D}$
corresponding to the eigenvalues $1, -1$, respectively. Then $M$ is
normal
 if and only if $T^+$ and $T^-$ are involutive and
$\xi$ is foliate with respect to both $T^{+}$ and $T^{-}$.
\end{proposition}
\begin{proof}
Assume that $N^{(1)}$ vanishes identically. In particular, for any
$X,Y\in\Gamma(T^{+})$ we have
\begin{align*}
0&=[X,Y]-\eta([X,Y])\xi+[X,Y]-\psi[X,Y]-\psi[X,Y]-2d\eta(X,Y)\xi\\
&=2[X,Y]-2\psi[X,Y],
\end{align*}
so that the integrability of $T^{+}$ follows. Moreover, for any
$X\in\Gamma(T^{+})$
\begin{equation*}
0=N^{(1)}(X,\xi)=\psi^2[X,\xi]-\psi[\psi
X,\xi]=-\psi[X,\xi]+[X,\xi],
\end{equation*}
and hence $[X,\xi]\in\Gamma(T^{+})$. In a similar way the analogous
assertions for $T^{-}$ can be proved. Conversely, assume that $T^+$
and $T^-$ are both integrable. Then for all $X,Y\in\Gamma(T^{+})$ we
have
\begin{align*}
N^{(1)}(X,Y)&=\psi^2[X,Y]+[\psi X,\psi Y]-\psi[\psi X,Y]-\psi[X,\psi
Y]-2d\eta(X,Y)\xi\\
&=[X,Y]+[X,Y]-[X,Y]-[X,Y]=0
\end{align*}
and, analogously, $N^{(1)}(X,Y)=0$ for all $X,Y\in\Gamma(T^{-})$.
Next, if $X\in\Gamma(T^{+})$ and $Y\in\Gamma(T^{-})$,
\begin{align*}
N^{(1)}(X,Y)&=\psi^2[X,Y]+[\psi X,\psi Y]-\psi[\psi X,Y]-\psi[X,\psi
Y]-2d\eta(X,Y)\xi\\
&=[X,Y]-\eta([X,Y])\xi-[X,Y]-\psi[X,Y]+\psi[X,Y]-2d\eta(X,Y)\xi=0.
\end{align*}
Finally, for all $X\in\Gamma(T^{+})$ we have
$N^{(1)}(X,\xi)=[X,\xi]-\psi[X,\xi]=0$, since $\xi$ is foliate with
respect to $T^{+}$. Similarly one has $N^{(1)}(Y,\xi)=0$ for all
$Y\in\Gamma(T^{-})$.
\end{proof}

In a paracontact metric manifold one defines a symmetric, trace-free
operator $h:=\frac{1}{2}{\mathcal L}_{\xi}\psi$. $h$ anti-commutes
with $\psi$ and satisfies $h\xi=0$ and $\nabla\xi=-\psi  + \psi h $.
Moreover $h\equiv 0$ if and only if $\xi$ is a Killing vector field
and in this case $(M,\psi,\xi,\eta,g)$ is said to be a
\emph{$K$-paracontact manifold}. A normal paracontact metric
manifold is called a \emph{para-Sasakian manifold}. We have the
following characterization.

\begin{theorem}[\cite{zamkovoy}]\label{parasasakian1}
An almost paracontact metric structure $(\psi,\xi,\eta,g)$ is
para-Sasakian if and only if, for all $X,Y\in\Gamma(TM)$,
\begin{equation}\label{condizioneparasasaki}
(\nabla_{X}\psi)Y=-g(X,Y)\xi+\eta(Y)X.
\end{equation}
In particular, any para-Sasakian manifold is $K$-paracontact.
\end{theorem}

In any paracontact metric manifold S. Zamkovoy introduced a
canonical connection which  plays in paracontact geometry the same
role  of the generalized Tanaka-Webster connection (\cite{tanno}) in
a contact metric manifold.

\begin{theorem}[\cite{zamkovoy}]\label{paratanaka}
On a paracontact metric manifold there exists a unique connection
${\nabla}^{pc}$, called the \emph{canonical paracontact connection},
satisfying the following properties:
\begin{enumerate}
  \item[(i)] ${\nabla}^{pc}\eta=0$, ${\nabla}^{pc}\xi=0$,
  ${\nabla}^{pc}g=0$,
  \item[(ii)]
  $(\nabla^{pc}_{X}\psi)Y=(\nabla_{X}\psi)Y+g(X-hX,Y)\xi-\eta(Y)(X-hX)$,
  \item[(iii)] $T^{pc}(\xi,\psi Y)=-\psi T^{pc}(\xi,Y)$,
  \item[(iv)] $T^{pc}(X,Y)=2d\eta(X,Y)\xi$ on ${\mathcal D}=\ker(\eta)$.
\end{enumerate}
The explicit expression of this connection is given by
\begin{equation}\label{paradefinition}
\nabla^{pc}_{X}Y=\nabla_{X}Y+\eta(X)\psi Y+\eta(Y)(\psi X-\psi h
X)+g(X-hX,\psi Y)\xi.
\end{equation}
Moreover, the tensor torsion field is given by
\begin{equation}\label{paratorsion}
T^{pc}(X,Y)=\eta(X)\psi h Y - \eta(Y)\psi h X + 2g(X,\psi Y)\xi.
\end{equation}
\end{theorem}

An almost paracontact structure $(\psi,\xi,\eta)$ is said to be
\emph{integrable} (\cite{zamkovoy}) if the almost paracomplex
structure $\psi|_{\mathcal D}$ satisfies the condition
$N_{\psi}(X,Y)\in\Gamma(\mathbb{R}\xi)$  for all
$X,Y\in\Gamma(\mathcal D)$. This is equivalent to require that the
eigendistributions $T^{\pm}$ of $\psi$ are formally integrable, in
the sense that $[T^{\pm},T^{\pm}]\subset
T^{\pm}\oplus\mathbb{R}\xi$. For an integrable paracontact metric
manifold, the canonical paracontact connection shares many of the
properties of the Tanaka-Webster connection on CR-manifolds. For
instance we have the following result.

\begin{theorem}[\cite{zamkovoy}]\label{integrability}
A paracontact metric manifold $(M,\psi,\xi,\eta,g)$ is integrable if
and only if the canonical paracontact connection preserves the
structure tensor $\psi$.
\end{theorem}

In particular, from Theorem \ref{integrability} and Theorem
\ref{parasasakian1} it follows that any para-Sasakian manifold is
integrable.

\section{The main results}

In this section we prove that paracontact geometry and the theory of
almost bi-Legendrian manifolds are strictly connected. More
precisely we have the following result.

\begin{theorem}\label{main1}
There is a biunivocal correspondence between almost bi-Legendrian
structures and paracontact metric structures.
\end{theorem}
\begin{proof}
Let $(M,\eta,L_1,L_2)$ be an almost bi-Legendrian manifold. We
define a $(1,1)$-tensor field $\psi$ on $M$ by setting
$\psi|_{L_1}=I$, $\psi|_{L_2}=-I$ and $\psi\xi=0$. Moreover we put
\begin{equation}\label{metricacanonica}
g(X,Y):=d\eta(X,\psi Y)+\eta(X)\eta(Y)
\end{equation}
for all $X,Y\in\Gamma(TM)$. We prove that $(M,\psi,\xi,\eta,g)$ is
in fact a paracontact metric manifold. A straightforward computation
shows that $\psi^2=I-\eta\otimes\xi$. Moreover, by construction,
$\psi$ induces an almost paracomplex structure on $\mathcal D$ since
$\psi|_{\mathcal D}^2=I$ and the eigendistributions corresponding to
the eigenvalues $1, -1$ of $\psi|_{\mathcal D}$ are  $L_1$ and
$L_2$, respectively. It remains to check that $g$ is a compatible
metric such that $g(X,\psi Y)=d\eta(X,Y)$. Indeed we have, for any
$X,Y\in\Gamma(TM)$,
\begin{align*}
g(\psi X,\psi Y)&=d\eta(\psi X,Y-\eta(Y)\xi)+\eta(\psi X)\eta(\psi
Y)\\
&=d\eta(\psi X,Y)\\
&=-g(X,Y)+\eta(X)\eta(Y).
\end{align*}
Then $g(X,\psi Y)=d\eta(X,\psi^2 Y)+\eta(X)\eta(\psi
Y)=d\eta(X,Y-\eta(Y)\xi)=d\eta(X,Y)$ for all $X,Y\in\Gamma(TM)$.
Conversely, let $(M,\psi,\xi,\eta,g)$ be a paracontact metric
manifold. Then $\eta\wedge(d\eta)^n\neq 0$ everywhere on $M$
(\cite{zamkovoy}), that is $\eta$ is a contact form on $M$. We
define an almost bi-Legendrian structure $(L_1,L_2)$ on $M$ by
setting $L_1:=T^{+}$ and $L_2:=T^{-}$. Note that $T^{+}$ and $T^{-}$
are $g$-isotropic distributions, that is $g(X,Y)=0$ for all
$X,Y\in\Gamma(T^{\pm})$. Indeed by \eqref{compatibile} one has
$g(X,Y)=g(\psi X,\psi Y)=-g(X,Y)$, so that $g(X,Y)=0$. Consequently
$L_1$ and $L_2$ are in fact Legendrian distributions on $M$, since
they are $n$-dimensional subbundles of $\mathcal D=\ker(\eta)$
satisfying $d\eta(X,X')=0$  for all $X,X'\in\Gamma(L_1)$,
$d\eta(Y,Y')=0$  for all $Y,Y'\in\Gamma(L_2)$, and they are mutually
transversal.
\end{proof}

From now on, we shall identify a paracontact metric manifold  with
the canonical induced almost bi-Legendrian structure, according to
Theorem \ref{main1}. In the next results, we shall investigate how
that bijection acts on special classes of paracontact metric
manifolds.

\begin{corollary}\label{Kparacontact}
There \ is \ a \ biunivocal \ correspondence \ between \ flat \
almost bi-Legendrian structures and $K$-paracontact structures.
\end{corollary}
\begin{proof}
First, we give an explicit expression of the tensor field $h$ in
terms of the bi-Legendrian structure $(L_1,L_2)$ which is identified
with the pair of eigendistributions $(T^{+},T^{-})$ of the
corresponding paracontact metric structure $(\psi,\xi,\eta,g)$,
according to Theorem \ref{main1}. Indeed for any $X\in\Gamma(L_1)$,
we have
\begin{align}
h X&=\frac{1}{2}([\xi,X]-\psi([\xi,X]_{L_1}+[\xi,X]_{L_2}) \nonumber \\
&=\frac{1}{2}([\xi,X]-[\xi,X]_{L_1}+[\xi,X]_{L_2})\label{formulah}\\
&=[\xi,X]_{L_2}.\nonumber
\end{align}
Analogously  one obtains that $h Y = -[\xi, Y]_{L_1}$ for all
$Y\in\Gamma(L_2)$. Therefore $(\psi,\xi,\eta,g)$ is a
$K$-paracontact  structure if and only if both $L_1$ and $L_2$ are
flat as Legendrian distributions. This remark together with Theorem
\ref{main1} proves the assertion.
\end{proof}

\begin{corollary}\label{integ}
There is a biunivocal correspondence between bi-Legendrian
structures and integrable paracontact metric structures.
\end{corollary}
\begin{proof}
We have only to prove that an integrable paracontact metric
structure gives rise to a bi-Legendrian structure. Let
$(M,\psi,\xi,\eta,g)$ be an integrable paracontact metric manifold.
According to Theorem \ref{main1} we define an almost bi-Legendrian
structure on $(M,\eta)$ by setting $L_1:=T^{+}$ and $L_2:=T^{-}$.
Since $(\psi,\xi,\eta)$ is integrable we have $[L_{i},L_{i}]\subset
L_{i}\oplus\mathbb{R}\xi$, $i\in\left\{1,2\right\}$. But $L_i$ are
Legendrian distributions, so that $[L_i,L_i]\subset\mathcal D$,
$i\in\left\{1,2\right\}$.  Hence $L_1$ and $L_2$ are necessarily
involutive.
\end{proof}

Finally, the following result, which is a direct consequence of
Theorem \ref{main1} and Proposition \ref{parasasakian}, holds.

\begin{corollary}\label{parasasaki}
There is a biunivocal correspondence between flat bi-Legendrian
structures and para-Sasakian structures.
\end{corollary}

We recall that, given an almost bi-Legendrian manifold
$(M,\eta,L_1,L_2)$, an \emph{almost bi-Legendrian equivalence} is a
contactomorphism $f$ of $M$ which preserves the Lagrangian
distributions $L_1$ and $L_2$, that is $f^{\ast}\eta=\eta$,
${f_\ast}_x\xi_x=\xi_{f(x)}$ and ${f_\ast}_x({L_1}_x)={L_1}_{f(x)}$,
${f_\ast}_x({L_2}_x)={L_2}_{f(x)}$ for all $x\in M$ (\cite{pang}).
However, in view of Theorem \ref{main1} this notion reflects in the
corresponding paracontact structure as follows.

\begin{proposition}\label{equivalenza}
Let $(M,\eta,L_1,L_2)$ be an almost bi-Legendrian manifold. Then a
diffeomorphism $f:M\longrightarrow M$ is an almost bi-Legendrian
equivalence if and only if it is an automorphism with respect to the
induced paracontact metric structure.
\end{proposition}
\begin{proof}
\ We \ recall \ that \ an \ automorphism \  of \ a \ paracontact \
metric \ manifold \ $(M,\psi,\xi,\eta,g)$ is nothing but an isometry
$f:M\longrightarrow M$ such that $\psi\circ f_\ast=f_\ast\circ\psi$
and $f^\ast\eta=\eta$. Now let $f$ be an almost bi-Legendrian
equivalence. Then, by the identification $L_1\equiv T^{+}$,
$L_2\equiv T^{-}$ one has, for all $X\in\Gamma(L_1)$ and
$Y\in\Gamma(L_2)$, $\psi {f_\ast}X={f_\ast}X={f_\ast}\psi X$ and
$\psi {f_\ast}Y=-{f_\ast}Y={f_\ast}\psi Y$, so that, since
$f_\ast\xi=\xi$, it follows that $\psi\circ f_\ast=f_\ast\circ\psi$.
Then, for all $X,Y\in\Gamma(TM)$, $g(f_\ast X,f_\ast Y)=d\eta(f_\ast
X,\psi f_\ast Y)=d\eta(f_\ast X,f_\ast\psi Y)=(f^\ast d\eta)(X,\psi
Y)=d\eta(X,\psi Y)=g(X,Y)$. Conversely, if $f$ is an automorphism of
the paracontact metric manifold $(M,\psi,\xi,\eta,g)$, then for any
$X\in\Gamma(L_1)$ we have $\psi {f_\ast}X={f_\ast}\psi X={f_\ast}X$
and, consequently, $f_\ast X \in\Gamma(L_1)$. Analogously one can
prove that $f$ preserves the Legendrian distribution $L_2$.
\end{proof}

In particular, Theorem \ref{main1} and Proposition \ref{equivalenza}
permits to provide a classification of paracontact metric structures
in the following way. Let $(\psi,\xi,\eta,g)$ be a paracontact
metric structure on the contact manifold $(M,\eta)$. By Theorem
\ref{main1} we can associate with $(\psi,\xi,\eta,g)$ a canonical
almost bi-Legendrian structure $(L_1,L_2)$. Hence we can consider
the bilinear forms $\Pi_{L_1}$ and $\Pi_{L_2}$ of each Legendrian
distribution, according to $\S$ \ref{prefoliations}. The bilinear
forms $\Pi_{L_1}$ and $\Pi_{L_2}$ are invariant under almost
bi-Legendrian equivalences and hence, by Proposition
\ref{equivalenza}, also under automorphisms of the paracontact
metric structure $(\psi,\xi,\eta,g)$. Therefore $\Pi_{L_1}$ and
$\Pi_{L_2}$ may be considered as invariants of the given paracontact
metric structure $(\psi,\xi,\eta,g)$ and they can be easily used for
classifying paracontact metric structures according to the flatness,
degeneracy or non-degeneracy of $\Pi_{L_1}$ and $\Pi_{L_2}$. A
further criterion is the integrability/non-integrability of each of
the Legendrian distributions $L_1$ and $L_2$. For instance, we have
the class given by those paracontact metric structures such that the
induced almost bi-Legendrian structure is integrable and flat. By
Corollary \ref{parasasaki} this class corresponds to para-Sasakian
structures.   The total number of classes amounts to $36$.

\medskip

The study of the interplays between paracontact and bi-Legendrian
manifolds is also motivated by the following theorem, whose
consequences will be discussed  in the sequel.

\begin{theorem}\label{connection}
\ Let \ $(M,\eta,L_1,L_2)$ \ be \ an \ almost \ bi-Legendrian \
manifold \ and  \ let $(\psi,\xi,\eta,g)$ be the canonical
paracontact metric structure induced on $M$, according to Theorem
\ref{main1}. Let $\nabla^{bl}$ and $\nabla^{pc}$ be the
corresponding bi-Legendrian and canonical paracontact connections.
Then
\begin{enumerate}
  \item[(a)] $\nabla^{bl}\psi=0$, $\nabla^{bl}g=0$;
  \item[(b)] the bi-Legendrian and the canonical paracontact connections coincide if and only if
  the induced paracontact metric structure is integrable.
\end{enumerate}
\end{theorem}
\begin{proof}
(a) \ For any $V\in\Gamma(TM)$, $X\in\Gamma(L_1)$, $Y\in\Gamma(L_2)$
one has
\begin{equation*}
(\nabla^{bl}_{V}\psi)X=\nabla^{bl}_{V}\psi
X-\psi\nabla^{bl}_{V}X=\nabla^{bl}_{V}X-\nabla^{bl}_{V}X=0
\end{equation*}
and
\begin{equation*}
(\nabla^{bl}_{V}\psi)Y=\nabla^{bl}_{V}\psi
Y-\psi\nabla^{bl}_{V}Y=-\nabla^{bl}_{V}Y+\nabla^{bl}_{V}Y=0,
\end{equation*}
since $\nabla^{bl}$ preserves $L_1$ and $L_2$. Moreover,
$(\nabla^{bl}_{V}\psi)\xi=\nabla^{bl}_{V}\psi\xi-\psi\nabla^{bl}_{V}\xi=0$,
since $\nabla^{BC}\xi=0$. Then, from the fact that
$\nabla^{bl}\psi=0$, $\nabla^{bl}\eta=0$ and $\nabla^{bl}d\eta=0$,
we have for all $X,Y,Z\in\Gamma(TM)$
\begin{align*}
(\nabla^{bl}_{X}g)(Y,Z)&=X(d\eta(Y,\psi Z)+\eta(Y)\eta(Z))-d\eta(\nabla^{bl}_{X}Y,\psi Z)-d\eta(Y,\psi\nabla^{bl}_{X}Z)\\
&=(\nabla^{bl}_{X}d\eta)(Y,\psi
Z)+(\nabla^{bl}_{X}\eta)(Y)\eta(Z)+\eta(Y)(\nabla^{bl}_{X}\eta)(Z)=0.
\end{align*}
(b)  \ We have to check that the bi-Legendrian connection of
$(M,\eta,L_1,L_2)$ satisfies the statements (i)--(v) of Theorem
\ref{paratanaka}. By Proposition \ref{proprieta} and the definition
of $\nabla^{bl}$ we have $\nabla^{bl}\xi=\nabla^{bl}\eta=0$ and, by
(a), $\nabla^{bl}g=0$. Taking the properties of the torsion tensor
field, (iii) of Theorem \ref{biconnection} and (iv) of Theorem
\ref{paratanaka} into account, we have  that
$T^{bl}(X,Y)=2d\eta(X,Y)\xi=T^{pc}(X,Y)$ for all $X\in\Gamma(L_1)$
and $Y\in\Gamma(L_2)$, and,  by \eqref{paratorsion},
\eqref{formulah},
\begin{gather*}
T^{pc}(X,\xi)=-\psi h X=h\psi X=h X=[\xi,X]_{L_2}=T^{bl}(X,\xi),\\
T^{pc}(Y,\xi)=-\psi h Y=h\psi Y=-h Y=[\xi,Y]_{L_1}=T^{bl}(Y,\xi).
\end{gather*}
Next, by Proposition \ref{proprieta} we have that, for any
$X,X'\in\Gamma(L_1)$, $Y,Y'\in\Gamma(L_2)$, $T^{bl}(X,X')$ \ $=$ \
$-[X,X']_{L_2}$, \ $T^{bl}(Y,Y')$ \ $=$ \ $-[Y,Y']_{L_1}$, \ and \
$T^{pc}(X,X')$ \ $=$ \ $2d\eta(X,X')\xi=0$, \
$T^{pc}(Y,Y')=2d\eta(Y,Y')\xi=0$, so that the torsion tensor fields
of the two connections $\nabla^{bl}$ and $\nabla^{pc}$ coincide if
and only if the Legendrian distributions $L_1$ and $L_2$ are
involutive. In view of Corollary \ref{integ} this is equivalent to
the integrability of the induced paracontact metric structure.
Finally, by virtue of (a), the bi-Legendrian connection satisfies
(ii) of Theorem \ref{paratanaka} if and only if $\nabla^{pc}\psi=0$.
But, by Theorem \ref{integrability}, that last condition is
equivalent to the integrability of the paracontact metric structure.
The theorem is thus completely proved.
\end{proof}

\begin{corollary}\label{clarification}
\ Any bi-Legendrian manifold $(M,\eta,{\mathcal F}_1,{\mathcal
F}_2)$ admits a canonical paracontact metric structure whose
paracontact connection coincides with the bi-Legendrian connection
of $(M,\eta,{\mathcal F}_1,{\mathcal F}_2)$.
\end{corollary}

\begin{remark}
Note that the explicit definition of the bi-Legendrian connection,
which we have recalled in $\S$ \ref{prefoliations}, is rather
involved. Thus Corollary \ref{clarification} clarifies its nature:
the bi-Legendrian connection can be seen as the canonical
paracontact connection of an integrable paracontact metric
structure. Moreover, by \eqref{paradefinition}, we can also deduce a
 relation between the bi-Legendrian connection and the Levi-Civita
connection.
\end{remark}

\medskip

Theorem \ref{main1} and its  corollaries yield some consequences for
the theory of paracontact manifolds, especially for para-Sasakian
geometry. Indeed many of the known results about bi-Legendrian
manifolds may be transported to the induced paracontact metric
structure. A crucial role in such matters is played by the relations
between the bi-Legendrian and the canonical paracontact connection
proved in Theorem \ref{connection}.

We recall that the bi-Legendrian connection of an almost
bi-Legendrian manifold $(M,\eta,L_1,L_2)$ is said to be
\emph{tangential} (\cite{mino2}) if its curvature tensor field
satisfies $R^{bl}(X,Y)=0$ for all $X\in\Gamma(L_1)$ and
$Y\in\Gamma(L_2)$. The geometric meaning of tangentiality is
explained in \cite{mino2}, where  some strong consequences  on the
geometry of the manifold are proved. Then we have the following
results.

\begin{proposition}
Let $(M,\psi,\xi,\eta,g)$ be a $(2n+1)$-dimensional K-paracontact
manifold. Suppose that one among $T^{+}$ and $T^{-}$ is integrable
and the associated bi-Legendrian connection is tangential. Then
$\emph{Pont}^j(TM)$ vanishes for $j>n$, where $\emph{Pont}(TM)$
denotes the Pontryagin algebra of the bundle $TM$.
\end{proposition}
\begin{proof}
  The \ assertion \ follows \ by \ applying \ \cite[Theorem 5.4]{mino2} \ to \ the \
 flat \ almost \ bi-Legendrian structure induced on $M$ by the $K$-paracontact
structure $(\psi,\xi,\eta,g)$ according to Corollary
\ref{Kparacontact}.
\end{proof}

\begin{proposition}\label{ehrconnection}
Let $(M,\psi,\xi,\eta,g)$ be a compact connected K-paracontact
manifold. Suppose that  $T^{+}$ (respectively, $T^{-}$) is
integrable and the corresponding bi-Legendrian connection is
tangential. If the leaves of the foliation defined by $T^{+}$
(respectively, $T^{-}$) are complete affine manifolds (with respect
to the bi-Legendrian connection), then the subbundle
$T^{-}\oplus\mathbb{R}\xi$ (respectively,
$T^{+}\oplus\mathbb{R}\xi$) is an Ehresmann connection, in the sense
of \cite{blumenthal1}, for the foliation defined by $T^{+}$
(respectively, $T^{-}$).
\end{proposition}
\begin{proof}
The assertion follows from \cite[Theorem 6.2]{mino2}.
\end{proof}

\begin{corollary}
Let $(M,\psi,\xi,\eta,g)$ be a compact connected para-Sasakian
manifold. Assume that the canonical paracontact connection
$\nabla^{pc}$ is tangential and the leaves of the foliation defined
by $T^{+}$ (respectively, $T^{-}$) are complete affine manifolds
with respect to $\nabla^{pc}$. Then the subbundle
$D^-=T^{-}\oplus\mathbb{R}\xi$ (respectively,
$D^+=T^{+}\oplus\mathbb{R}\xi$) is an Ehresmann connection for the
foliation defined by $T^{+}$ (respectively, $T^{-}$). Furthermore,
the universal covers of any two leaves of $T^+$ (respectively,
$T^{-}$) are isomorphic and the universal cover $\tilde M$ of $M$ is
topologically a product $\tilde{\mathcal L}^{+}\times\tilde{\mathcal
D}^{-}$ (respectively, $\tilde{\mathcal L}^{-}\times\tilde{\mathcal
D}^{+}$), where $\tilde{\mathcal L}^{+}$ (respectively,
$\tilde{\mathcal L}^{-}$) is the universal cover of the leaves of
$T^{+}$ (respectively, $T^{-}$) and $\tilde{\mathcal D}^{-}$
(respectively, $\tilde{\mathcal L}^{+}$) is the universal cover of
the leaves of the foliation defined by $D^{-}$ (respectively,
$D^{+}$).
\end{corollary}
\begin{proof}
The proof follows from Proposition \ref{ehrconnection} and Corollary
6.7, 6.9 in  \cite{mino2}, taking into account that $T^{+}$ and
$T^{-}$ are involutive and flat and, by Theorem \ref{connection},
the bi-Legendrian and the canonical paracontact connections
coincide.
\end{proof}

\begin{proposition}\label{leavesflat}
Let $(M,\psi,\xi,\eta,g)$ be a para-Sasakian manifold. Then the
canonical paracontact connection is flat along the leaves of the
foliations defined by the eigendistributions $T^{+}$ and $T^{-}$.
\end{proposition}
\begin{proof}
Since $M$ is para-Sasakian, the paracontact structure
$(\psi,\xi,\eta)$ is normal, in particular integrable. Hence by
Theorem \ref{connection} the canonical paracontact connection
$\nabla^{pc}$ coincides with the bi-Legendrian connection
$\nabla^{bl}$ of the induced flat bi-Legendrian structure defined by
$(T^{+},T^{-})$, according to Corollary \ref{parasasaki}. Then by
applying Proposition \ref{proprieta}  we get the result.
\end{proof}

\begin{corollary}
Let $(M,\psi,\xi,\eta,g)$ be a para-Sasakian manifold. Then the
leaves of the foliations defined by the eigendistributions $T^{+}$
and $T^{-}$ and the leaves of the foliations defined by
$T^{+}\oplus\mathbb{R}\xi$ and $T^{-}\oplus\mathbb{R}\xi$  admit a
canonical (flat) affine structure.
\end{corollary}
\begin{proof}
The first part is a direct consequence of Proposition
\ref{leavesflat}. We prove that the leaves of the foliations defined
by $T^{+}\oplus\mathbb{R}\xi$ and $T^{-}\oplus\mathbb{R}\xi$  admit
a canonical (flat) affine structure. Indeed, by \eqref{paratorsion}
the torsion of the canonical paracontact connection satisfies
$T^{pc}(Z,Z')=2d\eta(Z,Z')=-\eta([Z,Z'])=0$ for all
$Z,Z'\in\Gamma(T^{\pm})$  and $T^{pc}(Z,\xi)=\pm h Z=0$ for all
$Z\in\Gamma(T^{\pm})$, since the distributions $T^{+}$ and $T^{-}$
are integrable and $h=0$, $M$ being para-Sasakian. Moreover by
Theorem \ref{connection} $\nabla^{pc}=\nabla^{bl}$, hence using
Proposition \ref{proprieta}, we have that $R^{pc}(Z,\xi)=0$ for all
$Z\in\Gamma(T^{\pm})$. Thus the connection induced by the canonical
paracontact connection on the leaves of the foliations defined by
$T^{+}\oplus\mathbb{R}\xi$ and $T^{-}\oplus\mathbb{R}\xi$ provides
the desired flat affine structure.
\end{proof}

\section{Examples and remarks}

In this section we shall present some wide classes of examples of
bi-Legendrian manifolds and thus, in turn, again by Theorem
\ref{main1}, of paracontact metric structures. Moreover we prove a
theorem of local equivalence for $\nabla^{pc}$-flat paracontact
manifolds. Let us begin with a method for constructing new almost
bi-Legendrian structures from a given Legendrian distribution.

\subsection{Conjugate Legendrian foliation}
Let $L$ be a Legendrian distribution on a contact manifold
$(M,\eta)$. We show that there exists at least one Legendrian
distribution on $M$ transversal to $L$. Indeed, since $(M,\eta)$ is
a contact manifold it admits a Riemannian metric $G$ and a
$(1,1)$-tensor field $\phi$ satisfying
\begin{equation}\label{metrica1}
\phi^2  = -I + \eta\otimes\xi,  G(\phi X,\phi
Y)=G(X,Y)-\eta(X)\eta(Y),  G(X,\phi Y)=d\eta(X,Y)
\end{equation}
for all $X,Y\in\Gamma(TM)$ (cf. \cite{blair0}). $G$ is in fact an
associated metric ($\S$ \ref{prefoliations}). Note that from
\eqref{metrica1} it follows that
\begin{equation}\label{metrica2}
G(X,Y)=-d\eta(X,\phi Y)+\eta(X)\eta(Y).
\end{equation}
Let $Q$ be the distribution defined by $Q:=\phi L={\mathcal D}\cap
L^\perp$. As a matter of fact $Q$  is a Legendrian distribution on
$M$ called the \emph{conjugate Legendrian distribution} of $L$ with
respect to the contact metric structure $(\phi,\xi,\eta,G)$
(\cite{jayne1}). Thus $(L,Q)$ defines a $G$-orthogonal almost
bi-Legendrian structure on $M$. It should be remarked that under the
assumption that $\xi$ is Killing, if $L$ is flat, degenerate or
non-degenerate then also its conjugate is flat, degenerate or
non-degenerate, respectively (\cite{jayne1}).

\begin{proposition}
Let $L$ be a Legendrian distribution on a contact manifold
$(M,\eta)$ and let $Q$ be its conjugate Legendrian distribution.
Then starting from $(L,Q)$ one can define infinitely many almost
bi-Legendrian structures on $(M,\eta)$.
\end{proposition}
\begin{proof}
Let $(\psi,\xi,\eta,g)$ be the paracontact metric structure
 associated with the almost bi-Legendrian structure
$(L,Q)$. Note that
\begin{equation}\label{anticomm}
\phi\circ\psi=-\psi\circ\phi.
\end{equation}
Indeed, for any $X\in\Gamma(L)=\Gamma(T^{+})$, $\phi\psi X=\phi
X=-\psi\phi X$ and, for any $Y\in\Gamma(Q)=\Gamma(T^{-})$, $\phi\psi
Y=-\phi Y=-\psi\phi Y$, since $\phi L = Q$ and $\phi Q = L$.
Finally, $\phi\psi\xi=0=-\psi\phi\xi$. Now, let
$\alpha,\beta\in\mathbb{R}$ such that $\alpha^2+\beta^2=1$ and let
$\psi_{\alpha,\beta}$ be the tensor field of type $(1,1)$ defined by
$\psi_{\alpha,\beta}=\alpha\psi+\beta\phi\circ\psi$. One has, for
any $X\in\Gamma(TM)$,
\begin{align*}
\psi^2_{\alpha,\beta}X&=\alpha^2\psi^2 X+\alpha\beta\psi\phi\psi
X+\beta\alpha\phi\psi^2 X+\beta^2\phi\psi\phi\psi X\\
&=\alpha^2 X - \alpha^2\eta(X)\xi-\alpha\beta\phi X+\alpha\beta\phi
X+\beta^2 X-\beta^2\eta(X)\xi\\
&=X-\eta(X)\xi.
\end{align*}
Thus $(\psi_{\alpha,\beta},\xi,\eta)$ is an almost paracontact
 structure on $M$ and the eigenspace distributions
$T^{+}_{\alpha,\beta}$ and $T^{-}_{\alpha,\beta}$ of
$\psi_{\alpha,\beta}|_{\mathcal D}$ define an almost bi-Legendrian
structure on $M$. Indeed, let $X$ and $X'$ be two sections of
$T^{+}_{\alpha,\beta}$. So $X=\psi_{\alpha,\beta}X=\alpha\psi
X+\beta\phi\psi X$ and $X'=\psi_{\alpha,\beta}X'=\alpha\psi
X'+\beta\phi\psi X'$. Then we have
\begin{align*}
d\eta(X,X')&=\alpha^2 d\eta(\psi X,\psi X')+\alpha\beta d\eta(\psi
X,\phi\psi X')+\beta\alpha d\eta(\phi\psi X,\psi X')\\
&\quad+\beta^2 d\eta(\phi\psi
X,\phi\psi X')\\
&=-\alpha^2 d\eta(X,X')+\alpha\beta d\eta(X,\phi X')+\alpha\beta
d\eta(\phi X,X')-\beta^2 d\eta(X,X')\\
&=-d\eta(X,X'),
\end{align*}
so that $d\eta(X,X')=0$. Analogously one can prove that
$d\eta(Y,Y')=0$ for all $Y,Y'\in\Gamma(T^{-}_{\alpha,\beta})$.
\end{proof}

An interesting case occurs when $(M,\phi,\xi,\eta,G)$ is a Sasakian
manifold. We have in fact the following result.

\begin{theorem}
Let $(M,\phi,\xi,\eta,G)$   be a Sasakian manifold endowed with a
flat Legendrian distribution $L$. Let $Q=\phi L$ be its conjugate
Legendrian distribution. Suppose that the generalized Tanaka-Webster
connection $\nabla^{TW}$ preserves the distribution $L$. Then we
have:
\begin{enumerate}
  \item[(i)] $L$ and $Q$ are integrable and $\nabla^{TW}$
  coincides with the bi-Legendrian connection $\nabla^{bl}$
  corresponding to the bi-Legendrian structure $(L,Q)$.
  \item[(ii)] The induced paracontact metric structure $(\psi,\xi,\eta,g)$ is
  para-Sasakian and the distributions $T^{+}_{\alpha,\beta}$ and
  $T^{-}_{\alpha,\beta}$ are integrable.
  \item[(iii)] The canonical paracontact connection of
  $(\psi,\xi,\eta,g)$ coincides with the generalized
  Tanaka-Webster connection of $(\phi,\xi,\eta,G)$.
\end{enumerate}
\end{theorem}
\begin{proof}
The first part of the statement has been proved in \cite{mino3}.
Next, that the paracontact metric structure $(\psi,\xi,\eta,g)$ is
para-Sasakian follows from Corollary \ref{parasasaki}, taking into
account that $L$ and $Q$ are integrable and flat. We prove that the
Legendrian distributions $T^{\pm}_{\alpha,\beta}$, corresponding to
the almost paracontact structures $(\psi_{\alpha,\beta},\xi,\eta)$
defined above are involutive. For any
$X,X'\in\Gamma(T^{+}_{\alpha,\beta})$ we have, by
\eqref{condizionesasaki} and \eqref{condizioneparasasaki},
\begin{align*}
\psi_{\alpha,\beta}[X,X']&=\alpha\psi\nabla_{X}X' +
\beta\phi\psi\nabla_{X}X' -
\alpha\psi\nabla_{X'}X - \beta\phi\psi\nabla_{X'}X\\
&=-\alpha(\nabla_{X}\psi)X'+\alpha\nabla_{X}\psi
X'-\beta\phi(\nabla_{X}\psi)X'+\beta\phi\nabla_{X}\psi
X'\\
&\quad+\alpha(\nabla_{X'}\psi)X-\alpha\nabla_{X'}\psi
X+\beta\phi(\nabla_{X'}\psi)X-\beta\phi\nabla_{X'}\psi X\\
&=\alpha g(X,X')\xi-\alpha\eta(X')X+\alpha\nabla_{X}\psi X'+\beta
g(X,X')\phi\xi-\beta\eta(X')\phi X\\
&\quad-\beta(\nabla_{X}\phi)\psi X'+\beta\nabla_{X}\phi\psi
X'-\alpha g(X,X')\xi+\alpha\eta(X)X'-\alpha\nabla_{X'}\psi X\\
&\quad-\beta g(X,X')\phi\xi+\beta\eta(X)\phi
X'+\beta(\nabla_{X'}\phi)\psi X-\beta\nabla_{X'}\phi\psi X\\
&=\nabla_{X}\psi_{\alpha,\beta}X'-\beta g(X,\psi
X')\xi+\beta\eta(\psi X')X-\nabla_{X'}\psi_{\alpha,\beta}X\\
&\quad+\beta
g(X',\psi X)\xi-\beta\eta(\psi X)X'\\
&=\nabla_{X}X'-\nabla_{X'}X-2\beta d\eta(X,X')\xi\\
&=[X,X']
\end{align*}
since $d\eta(X,X')=0$, $T^{+}_{\alpha,\beta}$ being a Legendrian
distribution. By a similar argument one can prove also the
integrability of $T^{-}_{\alpha,\beta}$.  Thus (ii) is proved.
Finally, since $(\psi,\xi,\eta,g)$ is para-Sasakian, by Theorem
\ref{connection} the canonical paracontact connection coincides with
the bi-Legendrian connection of $(M,\eta,L,Q)$, which in turn, by
(i), coincides with the generalized Tanaka-Webster connection of the
Sasakian structure $(\phi,\xi,\eta,G)$.
\end{proof}

\subsection{Unit cotangent bundles}
Another way for attaching to a given Legendrian distribution $L$  a
transversal Legendrian distribution is proposed, by an intrinsic
construction, in \cite{pang} and \cite{libermann} under the
assumption of the integrability and non-degeneracy of $L$. So let
$\mathcal F$ be a non-degenerate Legendrian foliation. One defines
an operator $S_{\mathcal F}:{\mathcal D}\longrightarrow{\mathcal D}$
by setting $S_{\mathcal F}=\frac{1}{2}(i_{\mathcal D}-{\mathcal
L}_{\xi}\lambda)$, where $\lambda$ is the tensor field of type
$(1,1)$ on $M$ such that $\Pi(\lambda Z,X)=d\eta(Z,X)$ for all
$Z\in\Gamma(TM)$ and $X\in\Gamma(T\mathcal F)$. Then the image of
$S_{\mathcal F}$ is a Legendrian distribution of $(M,\eta)$
transversal to $\mathcal F$. In particular, that construction
applies to unit cotangent bundles. Let $M$ be a $(n+1)$-dimensional
smooth manifold with local coordinates
$\left(x_1,\ldots,x_{n+1}\right)$. Then $q_i=x_i\circ\pi$ and $p_i$,
$i\in\left\{1,\ldots,n+1\right\}$, are coordinates on the cotangent
bundle $T^{\ast}M$, where $(p_1,\ldots,p_{n+1})$ are fiber
coordinates and $\pi:T^{\ast}M\longrightarrow M$ denotes the
projection map. The Liouville form on $T^{\ast}M$ is then defined in
coordinates by $\beta=\sum_{i=1}^{n+1}p_i dq_i$. Now let
$F:T^{\ast}M\longrightarrow[0,+\infty[$ be a function such that
$F(tv) = t F(v)$ for all $t\geq 0$ and $v\in T^{\ast}M$. Then the
set $S^{\ast}_{F}M=\left\{v\in T^{\ast}M | F(v)=1\right\}$ is a
$(2n+1)$-dimensional submanifold of $T^{\ast}M$ called \emph{unit
cotangent bundle}. The Liouville form $\beta$ on $T^{\ast}M$
pulls-back to a contact form $\eta_{F}$ on $S^{\ast}_{F}M$ and the
connected components of the fibers of the projection
$\pi:S^{\ast}_{F}M\longrightarrow M$ define a Legendrian foliation
${\mathcal F}_F$ on $S^{\ast}_{F}M$. The computation in coordinates
of the invariant ${\Pi}_{{\mathcal F}_F}$ yields that it is the
restriction to the tangent bundle of ${\mathcal F}_F$ of the
symmetric form
\begin{equation*}
{\Pi}_{{\mathcal F}_F}=\sum_{i,j=1}^{n+1}\frac{\partial^2
F}{\partial p_i
\partial p_j}dp_i\otimes dp_j.
\end{equation*}
So  if  the  Hessian  matrix  $\left\{\frac{\partial^2 F}{\partial
p_i
\partial p_j}\right\}$  is  not singular then ${\mathcal F}_F$
defines a non-degenerate Legendrian foliation on $S^{\ast}_{F}M$.
Thus with the notation above we set $Q_{F}:=\textrm{Im}(S_{{\mathcal
F}_F})$ defining in this way a canonical almost bi-Legendrian
structure on the unit cotangent bundle $S^{\ast}_{F}M$.  For
example, if $F$ is the norm defined by a Riemannian metric $g$ on
$M$, then $S^{\ast}_{F}M$ is the cotangent sphere bundle
$S^{\ast}_{g}M=T^{\ast}_{1}M$ on $M$ and the Legendrian foliation
${\mathcal F}_g$ on $T^{\ast}_{1}M$ is clearly non-degenerate. By a
result of Pang (\cite[Proposition 5.30]{pang}) any positive definite
Legendrian foliation is locally equivalent to one of the form
${\mathcal F}_F$ with $F$ a Finslerian metric on $M$. In particular
this implies the following theorem.
\begin{proposition}
Let $\mathcal F$ be a positive definite Legendrian foliation on a
contact manifold $(M,\eta)$. Then the  almost bi-Legendrian manifold
$(M,\eta,T{\mathcal F}, Q_{\mathcal F})$, where $Q_{\mathcal
F}=\textrm{Im}(S_{\mathcal F})$, is locally equivalent to the almost
bi-Legendrian manifold $(S^{\ast}_{\mathcal F}M,\eta_{F},{\mathcal
F}_F,Q_F)$ for some Finslerian metric $F$ on $M$.
\end{proposition}

\subsection{Standard bi-Legendrian structure on
$\mathbb{R}^{2n+1}$}

A standard example of paracontact metric manifold is the following.
Consider in $\mathbb{R}^{2n+1}$ with global coordinates
$x_1,\ldots,x_n,y_1,\ldots,y_n,z$, the (1,1)-tensor field $\psi$
represented by the matrix
\begin{equation*}
 \left(
  \begin{array}{ccc}
    -I_n & 0 & 0 \\
    0 & I_n & 0 \\
    -y_i & 0 & 0 \\
  \end{array}
\right)
\end{equation*}
and put $\eta=dz-\sum_{i=1}^{n}y_idx_i$ and
$\xi=\frac{\partial}{\partial z}$. A straightforward computation
shows that $(\psi,\xi,\eta)$ defines an almost paracontact structure
on $\mathbb{R}^{2n+1}$. Moreover, one can consider the
semi-Riemannian metric $g$ given by
\begin{equation*}
 \left(
  \begin{array}{ccc}
    y_i y_j & \frac{1}{2}\delta_{ij} & -y_i \\
    \frac{1}{2}\delta_{ij} & 0 & 0 \\
    -y_j & 0 & 1 \\
  \end{array}
\right)
\end{equation*}
and check that $g$ is a compatible metric and
$(\mathbb{R}^{2n+1},\psi,\xi,\eta,g)$ is in fact a para-Sasakian
manifold. The corresponding bi-Legendrian structure on
$(\mathbb{R}^{2n+1},\eta)$ is given by the integrable Legendrian
distributions spanned by $\left\{X_1,\ldots,X_n\right\}$ and by
$\left\{Y_1,\ldots,Y_n\right\}$, where, for each
$i\in\left\{1,\ldots,n\right\}$, $X_i:=\frac{\partial}{\partial
y_i}$ and $Y_i:=\frac{\partial}{\partial
x_i}+y_i\frac{\partial}{\partial z}$. It is called the
\emph{standard bi-Legendrian structure} on $\mathbb{R}^{2n+1}$ (cf.
\cite{mino1}). One can verify that the canonical paracontact
connection is flat. Now we prove that, in a certain sense, also the
converse holds.

\begin{proposition}
Let $(M,\psi,\xi,\eta,g)$ be a para-Sasakian manifold. If the
canonical paracontact connection is everywhere flat, then
$(M,\psi,\xi,\eta,g)$ is locally isomorphic to the standard
paracontact metric structure of $\mathbb{R}^{2n+1}$.
\end{proposition}
\begin{proof}
 Let \ us \ suppose \ that \ the \ canonical \ paracontact \ connection \
$\nabla^{pc}$  \ of   $(M,\psi,\xi,\eta,g)$ is flat. By Theorem
\ref{connection} $\nabla^{pc}$ coincides with the bi-Legendrian
connection associated with the flat bi-Legendrian structure defined
by $(T^{+},T^{-})$. Thus, by applying \cite[Theorem 4.2]{mino1}, we
have that this bi-Legendrian structure is locally equivalent to the
standard bi-Legendrian structure on $\mathbb{R}^{2n+1}$. Now the
assertion follows directly from Proposition \ref{equivalenza}.
\end{proof}

\subsection{Anosov flows}

We recall for the convenience of the reader some definitions. Let
$M$ be a compact differentiable manifold. The flow $\{\omega_t\}$ of
a non-vanishing vector field $\xi$ on $M$ is said to be an
\emph{Anosov flow} (or $\xi$ to be an \emph{Anosov vector field}) if
there exist subbundles $E^s$ and $E^u$ which are invariant along the
flow and such that $TM=E^s\oplus E^u\oplus\mathbb{R}\xi$ and there
exists a Riemannian metric such that
\begin{gather}\label{condizioneanosov}
|{\omega_t}_{\ast}v|\leq a\exp(-ct)|v| \ \emph{ for all } t\geq 0
\emph{ and } v\in E^{s}_{x},\\
|{\omega_t}_{\ast}w|\leq a\exp(ct)|w| \ \emph{ for all } t\leq 0
\emph{ and } w\in E^{u}_{x},\nonumber
\end{gather}
where $a,c>0$ are constants independent of $x\in M$ and $v\in
E^{s}_{x}$, $w\in E^{u}_{x}$. The subbundles $E^s$ and $E^u$ are
called the \emph{stable} and \emph{unstable} subbundles. D. V.
Anosov proved that they are integrable and that also
$E^s\oplus\mathbb{R}\xi$ and $E^u\oplus\mathbb{R}\xi$ are integrable
(\cite{anosov}). Now let $(M,\phi,\xi,\eta,G)$ be a contact metric
manifold for which the Reeb vector field is Anosov. Then $(E^s,E^u)$
defines a bi-Legendrian structure on $M$. Moreover, the invariance
of $E^s$ and $E^u$ with respect to the flow can be expressed in
terms of Legendrian foliations just by the flatness of $(E^s,E^u)$.
Note that the paracontact metric structure induced on $M$ is in fact
para-Sasakian. Hence $\xi$ is Killing with respect to the
semi-Riemannian metric of the canonical paracontact metric structure
induced on $M$ by that bi-Legendrian structure, even if it can never
be Killing with respect to the associated metric $G$ in view of
\eqref{condizioneanosov}. The most notable example of a contact
manifold for which the Reeb vector field is Anosov is the tangent
sphere bundle of a negatively curved manifold.

\subsection{Contact metric $(\kappa,\mu)$-spaces}

Let $(M,\phi,\xi,\eta,G)$ be a contact metric $(\kappa,\mu)$-space,
that is a contact metric manifold satisfying
\begin{equation*}
R(X,Y)\xi=\kappa(\eta(Y)X-\eta(X)Y)+\mu(\eta(Y)\bar{h}X-\eta(X)\bar{h}Y)
\end{equation*}
for some constants $\kappa,\mu\in\mathbb{R}$, where $2\bar h$
denotes the Lie derivative of $\phi$ in the direction of $\xi$.
These manifolds have been introduced and deeply studied by D. E.
Blair, T. Kouforgiorgos and B. J. Papantoniou in \cite{blair1}. The
authors proved that necessarily $\kappa\leq 1$ and $\kappa=1$ if and
only if $M$ is Sasakian. Then for $\kappa<1$ the structure is not
Sasakian and $M$ admits three mutually orthogonal integrable
distributions ${\mathcal D}(0)=\mathbb{R}\xi$, ${\mathcal
D}(\lambda)$ and ${\mathcal D}(-\lambda)$ corresponding to the
eigenspaces of $\bar h$, where $\lambda=\sqrt{1-\kappa}$. Therefore
$(M,\eta,{\mathcal D}(\lambda),{\mathcal D}(-\lambda))$ is a
bi-Legendrian manifold. The induced paracontact metric structure is
integrable and thus the canonical paracontact and the bi-Legendrian
connections coincide. Moreover by a result in \cite{mino4} those
connections parallelize also $\phi$, $\bar h$ and $G$. The standard
example of a contact metric $(\kappa,\mu)$-space is given by the
tangent sphere bundle of a Riemannian manifold of constant sectional
curvature.\\
\\
\textbf{Acknowledgement} \ \ The author is very grateful to Prof.
Anna Maria Pastore for her kind suggestions and remarks which have
improved the article.


\small

\begin{thebibliography}{99}
\bibitem{alekseevski1}D. V. Alekseevski, V. Cort\'{e}s, A. S.
Galaev, T. Leistner, \textit{Cones over pseudo-Riemannian manifolds
and their holonomy}, J. Reine Angew. Math., to appear.

\bibitem{alekseevski2}D. V. Alekseevski, C. Medori, A. Tomassini,
\textit{Maximally homogeneous para-CR manifolds}, Ann. Glob. Anal.
Geom. \textbf{30} (2006), 1--27.

\bibitem{anosov}D. V. Anosov, \textit{Geodesic flows on closed Riemannian manifolds with negative
curvature}, Proc. Stlekov Inst. Math. \textbf{90} (1967).

\bibitem{blair1}D. E. Blair, T. Kouforgiorgos, B. J. Papantoniou,
\textit{Contact metric manifolds satisfying a nullity condition},
Israel J. Math. \textbf{91} (1995), 189--214.

\bibitem{blair0}D. E. Blair, \textit{Riemannian geometry of contact and symplectic
manifolds}, Birkh\"{a}user, Boston, 2002.

\bibitem{blumenthal1}R. Blumenthal, J. Hebda, \textit{Ehresmann connections for foliations}, Indiana Univ. Math. J. \textbf{33} n. 4 (1984), 597--611.

\bibitem{mino1}B. Cappelletti Montano, \textit{Bi-Legendrian connections},
Ann. Polon. Math. \textbf{86} (2005), 79--95.

\bibitem{mino2}B. Cappelletti Montano, \textit{Characteristic classes and Ehresmann connections for Legendrian
foliations}, Publ. Math. Debrecen \textbf{70} (2007), 395--425.

\bibitem{mino3}B. Cappelletti Montano, \textit{Some remarks on the generalized Tanaka-Webster connection of a contact metric
manifold}, Rocky Mountain J. Math., to appear.

\bibitem{mino4}B. Cappelletti Montano, L. Di Terlizzi,
\textit{Contact metric $(\kappa,\mu)$-spaces as bi-Legendrian
manifolds}, Bull. Austral. Math. Soc. \textbf{77} (2008), 373--386.

\bibitem{cortes1}V. Cort\'{e}s, C. Mayer, T. Mohaupt, F.
Saueressing, \textit{Special geometry of Euclidean supersimmetry 1.
Vector multiplets}, J. High Energy Phys. 03:028 (2004), 73 pp.

\bibitem{cortes2}V. Cort\'{e}s, M. A. Lawn, L. Sch\"{a}fer, \textit{Affine hyperspheres associated to special para-K\"{a}hler
manifolds}, Int. J. Geom. Methods Mod. Phys. \textbf{3} (2006),
995--1009.

\bibitem{olszak}P. Dacko, Z. Olszak, \textit{On weakly para-cosymplectic manifolds of dimension
3}, J. Geom. Phys. \textbf{57} (2007), 561--570.

\bibitem{erdem}S. Erdem, \textit{On almost (para)contact (hyperbolic) metric manifolds and harmonicity of $(\varphi,\varphi')$-holomorphic maps between
them}, Houston J. Math. \textbf{28} (2002), 21--45.

\bibitem{etayo1}F. Etayo, R. Santamaria, \textit{The canonical connection of a bi-Lagrangian
manifold}, J. Phys. A: Math. Gen. \textbf{34} (2001), 981--987.

\bibitem{zamkovoy2}S. Ivanov, D. Vassilev, S. Zamkovoy, \textit{Conformal paracontact curvature and the local flatness theorem}, preprint, available on arXiv:0707.3773v1 [math.DG].

\bibitem{jayne1}N. Jayne, \textit{Contact metric structures and Legendre
foliations}, New Zealand J. Math. \textbf{27} n. 1 (1998), 49--65.

\bibitem{kaneyuki1}S. Kaneyuki, F. L. Williams, \textit{Almost paracontact and parahodge structures on
manifolds}, Nagoya Math. J. \textbf{99} (1985), 173--187.


\bibitem{libermann}P. Libermann, \textit{Legendre foliations on
contact manifolds}, Different. Geom. Appl. \textbf{1} (1991),
57--76.

\bibitem{pang}M. Y. Pang, \textit{The structure of Legendre
foliations}, Trans. Amer. Math. Soc. \textbf{320} n. 2 (1990),
417--453.

\bibitem{tanno}S. Tanno, \textit{Variational problems on contact Riemannian
manifolds}, Trans. Amer. Math. Soc. \textbf{314} (1989), 349--379.

\bibitem{zamkovoy}S. Zamkovoy, \textit{Canonical connections on paracontact
manifolds}, preprint, available on arXiv:0707.1787v1 [math.DG].

\end{thebibliography}
\end{document}